\documentclass{math}

\english

\Author{Rastislav Telg\'arsky}

\Title{Dominant Frequency Extraction}

\ShortTitle{Dominant Frequency Extraction}

\Abstract{Time series are collected and studied extensively for the knowledge about the data source characteristics such as the trend or the spectral landscape. Some peaks in the spectral landscape correspond to dominant frequencies. The approach here is empirical: all time series are discrete and finite. Contents: Introduction. 1 Examples of periodic phenomena. 2 Algorithms and libraries. 3 Time series analysis. 4 Dominant frequency in ladar data. Conclusion. References. }


\Address{Department of Mathematics \\ Central New Mexico Community College \\ 525 Buena Vista Drive, Albuquerque, NM 87106, U.S.A.}{rtelgarsky@cnm.edu; rastislav@telgarsky.com}

\Subjclass {Primary 37M10, 62M15, 62P30, 37M05; Secondary 65T50.}  
 
\begin{document}

\SetUpTitle

\section*{Introduction}

While science tries to find the truth about the nature and the universe, the engineering approaches the properties of matter and phenomena which are necessary for human accomodation and survival, and which we can create and/or change. However, both science and engineering are simplifying our image of the world, sometimes into absurds. One way of simplification is to characterize a situation or a phenomenon by a single useful number, for example, a proportion, a volume, an average speed, an estimated probability, etc. Sometimes we find few numbers giving the useful characteristics. The observation of changes in the world leads to collections of data stamped by time; these data are called \textit{time series}. The time interval between data is mostly constant, but this restriction can be worked out.

We are interested here in phenomena called periodicity, periodic events, cyclical components, oscillations, vibrations, rhythms, resonance, seasonal variations, etc. The books on vibrations and waves which combine mathematics and nice illustrations are \cite{Crow,King}. There are too many papers on detection of periodic motion in video images; let's just mention these \cite{Cutler-Davis, Polana-Nelson, Seitz-Dyer}. There are several weaker notions of periodic functions which we do not consider here. These are: almost periodic, pseudo periodic, quasi periodic, nearly periodic, and semi periodic. The extraction of the dominant frequency is just a small part of knowledge extraction from time series \cite{Mandic}.

The temporal or dynamic changes in data can reflect a \textit{harmonic motion}, which has a sinusoid in some distorted form. A sine wave can be a projection of circular motion plot against time. A \textit{periodic motion} is all motion that repeats periodically. This includes the harmonic motion and pulses. While a pendulum of a wind-up clock performs a harmonic motion, the pendulum watched through a narrow window in front of the clock appears as pulses, which are periodic, but not harmonic. A \textit{random motion} occurs in erratic manner, as a sampling to an unknown probability distribution, and may contain all frequencies in a particular band. The main feature of a random motion is that it is not repeated or repeatable, in another words, there are no similarities of ups and downs (mountains and valleys) in different time intervals. The non-periodic non-random motions are piecewise linear or piecewise curved.

Here we consider (discrete) time series of finite length resulting from a data collection - measurement and recording. There are many situations when the observed data show a periodic behavior due to one frequency, called the \textit{dominant frequency}, which carries the maximum energy among all frequencies found in the spectrum. A similar notion is the \textit{fundamental frequency}, which is the smallest frequency having a peak among all frequencies in a power spectrum. The fundamental frequency is found in a vibrating string or organ pipe, along with weaker frequencies being the \textit{harmonics} (i.e., multiples) of the fundamental frequency (see \textit{harmonic} in \cite{Wiki}). The presence of harmonics in the signal demonstrates itself in the periodicity of the power spectrum, which has a comb of spikes. The usual technique for the detection of the fundamental frequency and its harmonics is the \textit{cepstrum}. There can be, however, other frequencies than harmonics present in the time series. They can be ordered by their energies, like it is done in electroencephalograms. We call them the \textit{2nd dominant}, the \textit{3rd dominant}, etc. When there are too many frequencies to characterize the periodicity of the data (time series), we can use the \textit{mean frequency}. The \textit{extraction of dominant frequency} means finding it, and sometimes removing it from the data as we remove linear trends. In that case, the plot of data with dominant frequency removed is called the \textit{residual plot}. It is sometimes more appropriate to replace the word "extraction" by \textit{determination} or \textit{estimation}. Other notions used for dominant frequency are dominating frequency, dominance frequency, dominant spike, or peak frequency. The dominant frequency can be the result of a \textit{resonance}, which can be desired (violins, vibraphones), or undesired (bridges, buildings), but it is always an important phenomenon. 

Whether the dominant frequency is seen in plots of data or not, the best way to reveal it is via a spectrogram (from Fourier transform) or a scalogram (form a wavelet transform). The \textit{spectral analysis} of data can reveal a hidden periodicity in data, so finding the dominant frequency is sometimes called \textit{spectrum peak picking}. The Fourier transform maps the data into the complex number domain. We can look at its real part and imaginary part, the amplitudes of frequencies and the phases of frequencies. The squares of amplitudes are called the \textit{power spectrum} or \textit{periodogram}. It is the computers and programming languages which started in 1950's the advanced data analysis a new science was born called \textit{Digital Signal Processing} (DSP), and little later, Digital Image Processing. 

\section{Examples of periodic phenomena} 
\begin{itemize}
\item A pendulum or any mechanical clock are suitable for testing some algorithms or demonstration of principles.
\item The heartbeat - the data is generated either by acoustic or electrical devices. The electrocardiogram (ECG) provides multiple time series. It is important for atrial fibrillation (see \cite{Wiki}) and ventricular fibrillation (see \cite{Wiki}) analyzed in cardiology.
\item The breathing - the activity of lungs is recorded by monitoring the flow of air, the movements of chest, or electrical impulses.
\item The electrical activities of brain are monitored by the electroencephalogram (EEG).
\item The electrical activities of digestive system are monitored by electrogastrogram (EGG).
\item The sleep laboratories have many different kinds of physiological devices and recorders, called polysomnograms, which include electroculograms (EOG) for recording of eye movements, and electromyograms (EMG) for the electrical activities of skeletal muscles.
\item The speech analysis has many uses. The analysis of acoustic recordings helps to recognize elements of speech for speech to text translations, and for logopedics.
\item The flapping of wings of birds and buzzing of insects. This is done with acoustic and/or video recording equipment.
\item The wavy motions of fish, octopus, crabs, and other sea creatures are studied in marine biology.
\item The elapsed-time photography helps to study slow-changing features in biology, geology, meteorology, etc.
\item The sunspots (protuberances) are monitored with multi-resolution cameras on the orbits of the Sun. We try to predict their occurrence and intensity, because of their negative impact on electromagnetic devices on Earth.
\item The vibrations of buildings, bridges, towers, constructions, machines, trucks, cars, etc. These are studied by vibrometry, which is both science and engineering domain.
\item The seasonal components of time series collected from industrial, commercial, financial, communication, transport of people and goods, energy consumption, road traffic at crossings, and other activities of human society.
\item The seismic activities of Earth and eruptions of volcanoes are not periodic when monitored in 24 hours/day, but they contain periodic segments of data, including infrasounds. The purpose of seismology and volcanology is to study the early signs of seismic activities in order to predict the earthquakes, and request evacuation of people in affected areas.
\item The meteor showers, the ocean waves, the ocean streams, and atmospheric winds are important for weather tracking. The weather patterns are too complicated to be characterized by time series. The observation of sunny and cloudy days is not sufficient for prediction of weather tomorrow and few days ahead. Even if they are reoccurring, they are not periodic events. However, the weather science (meteorology) makes constant progress with amount and speed of data processed for more accurate and useful forecasts.
\item The pattern recognition and texture characterization are branches of image processing, where the Fourier spectra are 2-dimensional, as well as other transforms and characteristics.
\end{itemize}

\section{Algorithms and libraries} 
In the empirical approach we use the professional software libraries in suitable programming languages to analyze the time series, and in particular, to detect the dominant frequencies.
The standard tools are: Fast Fourier Transform, Short-Time Fourier Transform, spectra, periodograms, and scalograms of various wavelets. These are included in many modern programming languages like C, C++, C\#, IDL, Java and Python. Then there are special languages for statistical and digital signal computing, like S-plus (a commercial implementation of the S programming language) and R language (free). There are many commercial software packages available for statistical analysis with great graphics. The most popular in research, engineering and education is Matlab from MathWorks \cite{MathWorks} with DSP, Image Processing, Wavelet, and other toolboxes. Moreover, there is constant flow of contributions from active users in all over the world. Good libraries and contributed codes are also available for Octave (free download), which is an alternative interpreter of Matlab m-files. For the time series occurring in economy there was developed the Berlin Procedure (Berlin Verfahren), which can perform seasonal adjustment to data, for example, daily, weekly, monthly or quarterly, and then perform an in-depth analysis. The latest version of the Berlin Procedure, BV4.1, is available for free download for non-commercial purposes. It turns out that market time series and not so random as they look in short time intervals. However, the prediction of time series is quite different problem than an extraction of the dominant frequency. Actually, the dominant frequencies in market data are call seasonal variations, such as monthly, quarterly and yearly. These are removed from the data, and then we have the residual time series, where we look for periodicities and other attributes.

\begin{figure}
  \begin{center}
    \includegraphics[width=10cm]{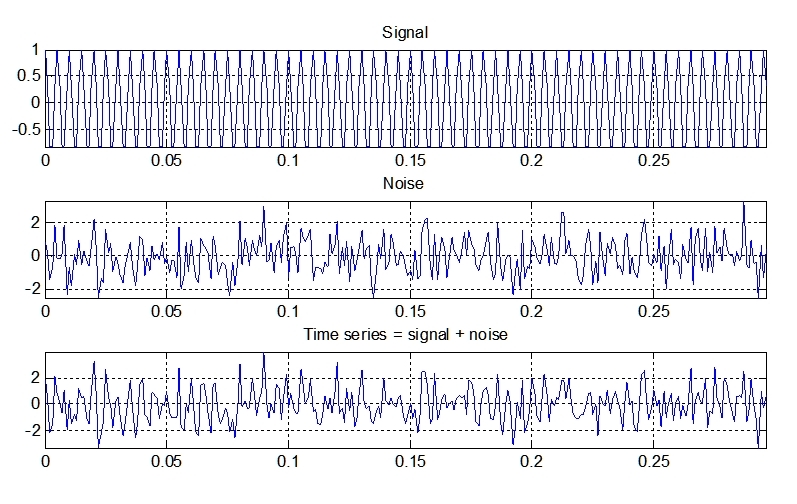}
  \end{center}
  \caption{The creation of time series with a dominant frequency.}
  \label{series}
\end{figure}

\begin{figure} 
  \begin{center}
    \includegraphics[width=10cm]{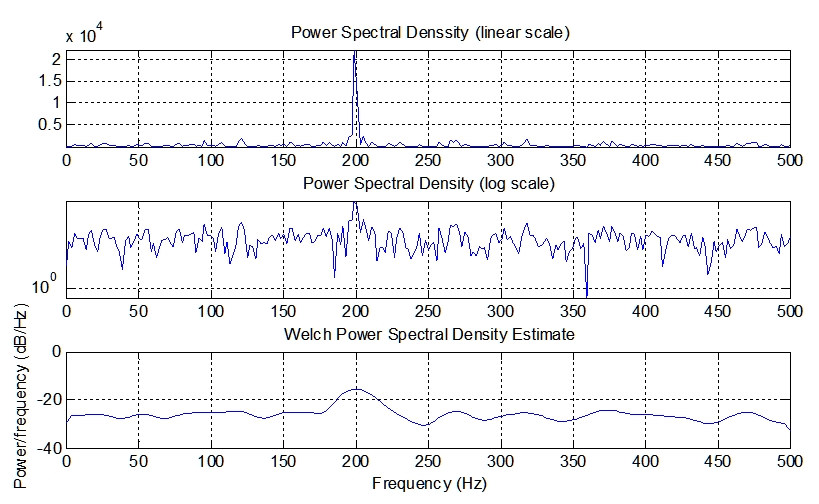}
  \end{center}
  \caption{The power spectra with the dominant frequency.}
  \label{psd}
\end{figure}

When a time series is collected from a relatively unknown source, it may have various components and their proportions what makes the extraction of dominant frequency quite difficult. On the other hand, we can create an artificial discrete time series of finite length by taking a finite combination of sinusoids representing different frequencies, and add a reasonably small Gaussian white noise. Then the sinusoid with the largest amplitude (the absolute value of the coefficient) is the dominating frequency. A typical entry level example consist of one sinusoid plus a Gaussian white noise (Figure~\ref{series}). The power spectrum of this time series has one pronounced peak clearly showing the dominant frequency (Figure~\ref{psd}). Here is the corresponding Matlab code:
\\
\\
\\
\begin{tt} 
\noindent{function demo\_dfe}\\
\noindent{Fs = 1000; \% sampling frequency 1 kHz}\\
\noindent{t = 0 : 1/Fs : 0.296; \% time scale}\\
\noindent{f = 200; \% Hz, embedded dominant frequency}\\
\noindent{x = cos(2*pi*f*t) + randn(size(t)); \% time series}\\
\noindent{plot(t,x), axis('tight'), grid('on'), title('Time series'), figure}\\
\noindent{nfft = 512; \% next larger power of 2}\\
\noindent{y = fft(x,nfft); \% Fast Fourier Transform}\\
\noindent{y = abs(y.\textasciicircum 2); \% raw power spectrum density}\\
\noindent{y = y(1:1+nfft/2); \% half-spectrum}\\
\noindent{[v,k] = max(y); \% find maximum}\\
\noindent{f\_scale = (0:nfft/2)* Fs/nfft; \% frequency scale}\\
\noindent{plot(f\_scale, y),axis('tight'),grid('on'),title('Dominant Frequency')}\\
\noindent{f\_est = f\_scale(k); \% dominant frequency estimate}\\
\noindent{fprintf('Dominant freq.: true \%f Hz, estimated \%f Hz\textbackslash n', f, f\_est)}\\
\noindent{fprintf('Frequency step (resolution) = \%f Hz\textbackslash n', f\_scale(2))}\\
\end{tt}

The thoughts about artificial mixtures of sinusoids and a white noise are not without a merit. The Pisarenko Harmonic Decomposition (PHD) does precisely that (see \cite{Wiki}), even more, it can tell amplitudes and frequencies of dominant frequencies up to certain number, which should be always specified ahead and should be small relative to the length of data. The Multiple Signal Classification (MUSIC) is a generalization of the Pisarenko Harmonic Decomposition. MUSIC accepts the complex-valued time series which are supposed to be the sum of $p$ complex exponentials and a complex Gaussian white noise. The algorithm returns $p$ largest peaks. The MUSIC algorithm is implemented in Matlab's Signal Processing Toolbox as \texttt{spectrum.music} \cite{Matlab-SP}.

\section{Time series analysis} 
The data (time series, signal) may come from another engineering team, or we might create them by running a possibly faithful simulation of the engineering process. Thus, in both cases, the data consist of finitely many samples, but they can be multi-dimensional, where the dimension is the number of features tracked. Let's focus on one dimensional case, where we study the temporal (dynamic) properties of the data source, and we focus on finding the dominating frequency. The task of dominant frequency extraction may look simple like this: find the peak in the power spectrum density of the signal to get the frequency, and then read the corresponding amplitude and phase from the Fourier transform of the signal. In cases where is a small number of fairly distant frequencies and a low value of the noise it is indeed so simple. However, in most important practical cases we need to discuss a number of issues.

\subsection{The characteristics}

$n_s$ = (number of samples in the data) = (length of the data sequence)\\
\noindent{$T_s$ = (time interval of collecting the samples) = (the observation time)}\\
\noindent{$r_s$ = (sampling rate) = $n_s / T_s$}\\
\noindent{$f_N$ = (Nyquist frequency) = $r_s / 2$}\\
\noindent{$f_{max}$ = (maximal frequency computable with FFT-based spectrum analysis) = $f_N$}\\
\noindent{$t_s$ = (time step) = (time between samples) = $T_s / n_s = 1 / r_s$}\\
\noindent{$p_f$  = (duration of one period of a frequency $f$) = $1 / f$ seconds}\\
\noindent{$np_f$ = (number of periods of a frequency $f$ during time $T_s$) = $T_s / p_f$ = $T_s  f$}\\
\noindent{$ns_f$ = (number of samples in one period of frequency $f$) = $n_s / np_f$}\\
\noindent{$f_{res}$ = (frequency resolution in Fourier spectrum analysis) = $1 / T_s$ Hz}\\
\noindent{$f_{min}$ = (minimal frequency computable via Fourier spectrum analysis) = $f_{res}$ Hz}

Alternatively, the sampling rate may be fixed, for example, $r_s$ = 30 Hz or 60 Hz, while the number of sample points $n_s$ is controlled by $T_s$: $n_s = r_s T_s$. If $n_s = 2 k$, where $k > 0$, then we can detect frequences $1/T_s, 2/T_s, ..., k/T_s$. 

\subsection{Example}

Let $n_s$ = 1000 samples and  $T_s$ = 100 seconds. Then

\noindent{$r_s$ = $n_s / T_s$ = 10 Hz.}\\
\noindent{$f_N$ = $r_s / 2$ = 5 Hz.}\\
\noindent{$f_{max}$ = 5 Hz.}\\
\noindent{$t_s$ = $1 / r_s$ = 0.1 second.}\\
Let f = 5 Hz. Then\\
\noindent{$p_f$  = $1 / f$ = 0.2 second.}\\
\noindent{$np_f$ = $T_s  f$ = 500 periods.}\\
\noindent{$ns_f$  = $n_s / np_f$ = 1000 / 500 = 2 samples / period.}\\
\noindent{$f_{res} = 1 / T_s$ = 0.01 Hz.}\\
\noindent{$f_{min}$ = 0.01 Hz.}

Note that the minimal frequency 0.01 Hz has the period of duration 100 seconds, and therefore it has only one period per entire duration of sampling. This frequency is supported by 1000 samples, and it is the smallest frequency computable (detectable) with Fourier transform methods. Therefore, the maximal resolution of its Fourier spectrum consists of frequencies: 0, 0.01, 0.02, ..., 5.0.

\subsection{Frequency limitations}

According to the Nyquist Theorem, we can sample and then estimate only those frequencies f which do not exceed $f_N$. In other words, we need at least 2 samples per period of a computable frequency. The minimal value of $n_s$ is 2, where the only detectable frequency is 1 Hz provided that $T_s = 1$ second. To have more samples is better, however, what we get is usually what we get, and there is no way to get more. The FFT methods also impose the limitation that $f \ge 1 / T_s$, or equivalently, we have for the maximum period $p_f \le T_s$. Therefore we must have at least one period per duration of the data collection. So, we get lower and upper bounds for $f$
\[
1 /T_s \le f \le (n_s / 2) / T_s = r_s / 2.
\]
In many applications, one period of a sinusoid does not indicate that this period will repeat many times to create a periodic event (phenomenon). I believe, as some other researchers do, that the sampling process must support at least 3 periods of the frequency. To detect smaller frequencies than $1 / T_s$ or even $3 / T_s$, we should either increase the observation time $T_s$ \emph{or} use different methods than Fourier spectral methods. Note that increasing $T_s$ decreases the Nyquist frequency when $n_s$ is fixed. 

\begin{figure}
  \begin{center}
    \includegraphics[width=8cm]{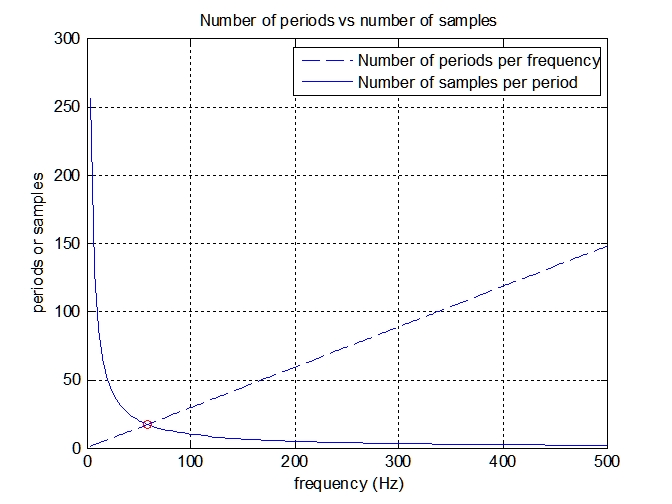}
  \end{center}
  \caption{Number of periods and number of samples.}
  \label{SDC}
\end{figure}

These limits are two extremes: either we have one period stretching over all time $T_s$ and containing all $n_s$ samples \emph{or} only only 2 samples per each period for all $n_s / 2$ periods. The number of periods of a frequency $f$ is $np_f = T_s  f$ and changes from 1 to $n_s / 2$ as $f$ increases. The number of samples per period of a frequency $f$ is $ns_f  = n_s / np_f = r_s / f$, and changes from $n_s$ down to 2 as $f$ increases. These two curves behave like the Supply and Demand Curves in the economic model of price determination in a market (see Figure~\ref{SDC}). It is easy to calculate that the equilibrium point is at the frequency $f = \sqrt{n_s} / T_s$, and this happens when $np_f  = ns_f = \sqrt{n_s}$.

However, the minimal frequency $f_{min}$, and the step in frequency scale of the spectrum $f_{res}$ are usually slightly different from $1 / T_s$. This  happens when the number of samples $n_s$ is not a power of 2. The processing of the number of points which are powers of 2 is required by the recursive formulas in the Fast Fourier Transform algorithm. Let $n_{fft}$ be the length of the Fast Fourier Transform used in the processing (this number is always a power of 2). For example, if $n_s = 296$, then  $n_{fft}$ can be chosen to be 256 or 512. Now, the formula for the minimal frequency and frequency resolution becomes
\[
f_{min} = f_{res} = r_s / n_{fft} = ( 1 / T_s) (n_s / n_{fft}).
\]
Clearly, if $n_s = n_{fft}$, then $f_{min} = 1 / T_s$.

\subsection{Signal components}

We try to deal with the signal as if it was made of 3 components: the trend, the waves (sinusoids) and the noise, leaving aside the component which is periodic, but with irregular periods. We think about trend as a polynomial, which is clearly a non-periodic component. The polynomial must be of a low degree, as is commonly used in non-linear (polynomial) regression. In the our applications, it was sufficient to deal with a linear trend. The Matlab function \texttt{detrend} removes the best straight-line fit linear trend from the data.
Each wave (sinusoid) of a digital signal representing time series is characterized by the frequency, the phase and the amplitude 
\[
     x(t) = A sin(2 \pi B t + C)
\]
The signal may contains several kind of noises, and not all of them additive. Noises may look like aperiodic multi-waves and their spectra may have many spikes. The effects of noises is usually decreased by filtering of the data. We can filter the signal with a \textit{band-pass filter} to cut out very low and very high frequency noises. However, the periodic part of the signal might also be weakened, because we may not know the band containing the dominant frequency. The Matlab function \textit{pwelch} performs the Power Spectral Density (PSD) estimate via Welch's method. By default of this algorithm, the time series is divided into eight sections with 50\% overlap, each section is windowed with a Hamming window, and eight modified periodograms are computed and averaged. This kind of spectrum is heavily doctored, so some people prefer the other extreme: the time series are either zero-padded to the nearest higher power of 2 or truncated to the nearest lower power of 2; then, FFT is taken and its the modulus is squared as in the Matlab demo code. This is the raw PSD. 

If some frequencies do not spread from the beginning of data to the end, then we must first determine the subintervals of their duration, because each dominant frequency must have specified the time of beginning and the time of end. This is typical, for example, for road vibrations, which are caused by passing of heavy trucks. In this case we use the Short-Time Fourier Transform (STFT) to determine the peaks and the subintervals, and then we can use the Fourier spectra for these subintervals.

\subsection{The peaks in PSD}

The next step of the algorithm is finding the local maxima in the power spectrum. The peaks (local maxima) in PSD can be detected with pixel accuracy by finding the indices of the array, or with a sub-pixel accuracy, for example, by fitting a smooth 'hat' over the peak in a small neighborhood. The peaks in PSD have approximately Gaussian curve shapes. The energy is stored not only in the dominant frequency but also in the width of the peak curve. Therefore, some sources recommend to consider a narrow bandwidth of dominant frequency, and define it as the width corresponding to frequencies with the amplitude $A/\sqrt 2$, where A is the amplitude of the peak, or use the Full Width at Half Maximum (FWHM) approach. For most applications, there are ready Matlab functions \texttt{findpeaks} and \texttt{localmax}.

However, in some cases finding the peaks is not good enough. What should be the relative power of the signal component with the dominant frequency? How much dominant should be the dominant frequency? If the entire power spectrum has the shape of a Gaussian distribution, then its peak is not necessarily strong enough to be considered as the dominant frequency. It may depend on the standard deviation of the spectrum; so one can apply some definitions of the signal to noise ratio (SNR). Some projects require that the peak corresponding to the dominant frequency is 30\% higher than other values in the power spectrum. One can ask how much energy (in \%) is contained in dominant frequency. Seitz and Dyer \cite{Seitz-Dyer} calculate the probability that a motion is periodic and compute the most likely period. Also, they use the Kolmogorov-Smirnov test to compare the near-periodic sequences. The Matlab function \texttt{spectrum.psd} can put a ribbon around the PSD estimate according to the specified confidence level $p$. Matlab has 3 functions to test the significance of the statement whether a signal is just a white noise or not: \texttt{kstest}, \texttt{lillietest} and \texttt{jbtest}.

Sometimes two or more peaks in PSD are too close, so that their separating requires a special attention. Matlab provides the useful function \texttt{periodogram} which may help to balance the filtering and sampling control parameters.

When plotting the phases of frequencies under the power spectrum, we can see, that due to the noise, phases may have jumps, in particular, the phase of the dominant frequency may have a jump. Thus, to extract the complete information from the power spectrum, we need to unwrap the phase array. For this phase correction there is a convenient Matlab function \texttt{unwrap}.

\subsection{Low frequencies}
If the dominating frequency is below $f_{min}$, then we need to apply different methods than spectral analysis. One approach is to use a least squares fit method to approximate the time series by the function 
\[
x(t) = A + B t + C \sin (2 \pi D t + E).
\]
I tried this when the sample data contained only 3 periods of the dominant frequency or less. When the time interval contained only a half of the period, the results were still good. When I reached below 1/4 of the period, the results were not so accurate due to the presence of white noise in the data. Matlab has the function \texttt{nlinfit}, and the Numerical Recipes \cite{NRC} have the steepest gradient algorithm \texttt{amoeba}.

All issues listed above indicate that the validation of extracted dominant frequency (with the amplitude and phase) must contain several factors such as the number of periods in the sample,  relative energy in the bandwidth of the dominant frequency, and the amount of white noise in the data. 

\section{Dominant frequency in ladar data} 
\begin{figure}
  \begin{center}
    \includegraphics[width=6cm]{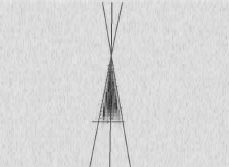}
  \end{center}
  \caption{Feature extraction from a range-Doppler ladar image.}
  \label{FE}
\end{figure}

The data to be analyzed for dominant frequencies were available only at the end of an engineering process \cite{Telg-et-al}, which I describe briefly. A moving and also rotating object was illuminated by an infrared laser, and its image was tracked by a special camera. This arrangement is called laser radar, or ladar. The laser was emitting trains of pulses and the camera recorded the time shifts between emitted and received pulses. Clearly, this was a preparation for the detection of the target range, and also for exploring the Doppler effect, that is, the shift in frequency caused by the component of the target motion in direction of the axis of illumination; actually, in all directions but those perpendicular to the axis of illumination. The resulting image had pixels with two coordinates: the time delay and the phase delay. This image was converted to another image, where the pixel coordinates were the distance and the velocity of the corresponding small portions of the target. This image is called the \textit{range-Doppler image}. The conversion was made by the Short-Time Fourier Transform \cite{Chen-Ling}. So, we got a video, where the target was tracked by a triangle or a quadrangle, and the locations of the vertices were recorded (see Figure~\ref{FE}). Now we got the multidimensional data, which we searched for dominant frequencies.

Not accidentally, ladar imitates the ultrasound imaging system of bats, which they use to catch bugs in fast flights. The only difference is in the variable resolution of images, because when a bat is closing the distance to its prey, which is flapping wings, the resolution is increasing to greater details.

The range-Doppler image is not quite similar to any standard optical image. The optical image has up-and-down and left-and-right orientations, while the range-Doppler image shows what is closer-and-farther and what is moving toward-and-away from the source of illumination. While ladars with their range-Doppler images are typically used for remote sensing applications, like satellites, the Doppler Radars have much wider applications \cite{Wiki}.

\begin{figure}
  \begin{center}
    \includegraphics[width=11cm]{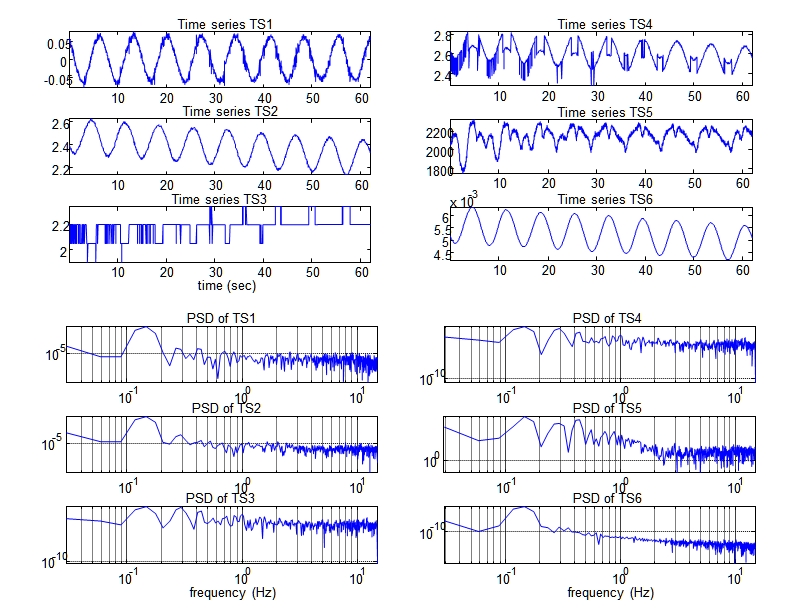}
  \end{center}
  \caption{Six time series and their power spectra.}
  \label{six-ts}
\end{figure}

An interesting part of the project was that I could use a computer simulation to generate a virtual target, where I could enter the dominating frequency as one of input parameters. However, in important test cases, the video data were created by another group, and I did not know the dominating frequency value prior to my analysis. My lucky number was that the dominating frequency, which I extracted from the time series, agreed to all specified decimal points with the frequency known only to the other group (Figure~\ref{six-ts}). 

There are many factors influencing the selection of processing methods for time series analysis. Some time series have to be processed in real time in order to trigger a prediction or even a warning, like the monitored vibrations of an overloaded bridge. There is a wide selection of Power Spectral Density estimators in Matlab \cite{Matlab-SP}: \texttt{pwelch, pmusic, pmtm, pcov, pmcov, pburg, pyulear, periodogram} and the function \texttt{spectrum}. The next group are wavelets \cite{Mallat}, where the time vs. scale resolution allows to detect events in certain time intervals and being significant in certain scale. The last group of algorithms I want to mention performs special decompositions of the signal such as the Empirical Mode Decomposition and the Independent Components Analysis (cf. \cite{Wiki}). The signal decomposition into multiple latent components using the latent Dirichlet allocation model is studied in \cite{Tel-Laf}.  

\section*{Conclusion}
The determination of dominant frequency helps to understand the structure of time series, and derive the consequences of the presence and intensity of the dominant frequency, whether it deals with biology, geology, medicine, or other fields of science and engineering. Recognizing the dominant frequency is a part of analysis of data leading to a better prediction, to more accurate diagnosis, and to a better tuned engineering design.

\end{document}